\documentclass[11pt,a4paper,reqno]{amsart}
\usepackage[utf8]{inputenc}
\usepackage[T1]{fontenc}
\usepackage{lmodern,microtype}
\usepackage{mathtools,amssymb,etoolbox,needspace}
\usepackage[margin=25mm]{geometry}
\usepackage[unicode,hidelinks]{hyperref}
\hypersetup{
 pdftitle={The equality $s(G)=\eta(G)+\exp(G)-1$ fails},
 pdfauthor={Guoqing Wang, Xiangneng Zeng},
 pdfsubject={The equality $s(G)=\eta(G)+\exp(G)-1$ fails},
 pdfkeywords={Zero-sum sequences, short zero-sum subsequences, Erd\H{o}s–Ginzburg–Ziv constant}
}
\allowdisplaybreaks[1]
\numberwithin{equation}{section}

\newtheorem{theorem}{Theorem}[section]
\newtheorem{proposition}[theorem]{Proposition}
\newtheorem{lemma}[theorem]{Lemma}

\theoremstyle{definition}

\theoremstyle{remark}

\newcommand{\F}{\mathbb F}

\newcommand{\sconst}{\mathsf{s}}
\newcommand{\Dconst}{\mathsf{D}}
\newcommand{\wt}{\operatorname{wt}}
\newcommand{\supp}{\operatorname{supp}}
\newcommand{\doi}[1]{\href{https://doi.org/#1}{\nolinkurl{#1}}}

\title[The equality $s(G)=\eta(G)+\exp(G)-1$ fails for $G=C_2^7\oplus C_6$]{The equality $s(G)=\eta(G)+\exp(G)-1$ fails for $G=C_2^7\oplus C_6$}
\author{Guoqing Wang}
\address{School of Mathematical Sciences \\ Tiangong University \\ Tianjin 300387, P.R. China}
\email{gqwang1979@aliyun.com}
\author{Xiangneng Zeng}
\address{Sino-French Institute of Nuclear Engineering and Technology,
Sun Yat-Sen University, Guangzhou 510275, P.R. China}
\email{zengxneng@mail.sysu.edu.cn}
\date{}
\makeatletter
\@ifundefined{subjclassname@2020}{%
 \@namedef{subjclassname@2020}{\textup{2020} Mathematics Subject Classification}}{}
\makeatother

\AtBeginEnvironment{theorem}{\samepage}
\AtBeginEnvironment{lemma}{\samepage}
\AtBeginEnvironment{proposition}{\samepage}

\begin{document}
\begin{abstract}
Let $G$ be a finite abelian group. The Erd\H{o}s–Ginzburg–Ziv constant $\sconst(G)$ is the smallest positive integer $\ell$ such that every sequence over $G$ of length at least $\ell$ contains a zero-sum subsequence of length $\exp(G)$. The short zero-sum constant $\eta(G)$ is defined similarly, with the requirement that the zero-sum subsequence is nonempty and has length at most $\exp(G)$.
It is known that
$\sconst(G)\ge \eta(G)+\exp(G)-1.$ A conjecture proposed by W. Gao \cite[Conjecture~2.3]{Gao03} in 2003, which has been widely studied in zero-sum theory, predicts that
$$
\sconst(G)=\eta(G)+\exp(G)-1
$$
for every nontrivial finite abelian group $G$. In this paper, we prove that this conjecture fails for the group $C_2^7\oplus C_6$.
\end{abstract}
\maketitle
\markboth{G. WANG, AND X. ZENG}{A COUNTEREXAMPLE TO GAO'S EQUALITY CONJECTURE}

\section{Introduction}
Let $G$ be a finite abelian group, and let $n=\exp(G)$.
A nonempty zero-sum subsequence is called \emph{short} if its length
is at most $n$.
The constant $\eta(G)$ is the least integer $\ell$ such that every
sequence over $G$ of length at least $\ell$ has a short zero-sum
subsequence.  The Erd\H{o}s--Ginzburg--Ziv constant $\sconst(G)$ is
defined in the same way, with the required length equal to $n$; Both invariants are widely studied in zero-sum theory, see~\cite[Section~6]{GG06} for a survey.

The celebrated Erd\H{o}s-Ginzburg-Ziv Theorem ~\cite{EGZ61} gives
$\sconst(C_n)=2n-1$; while $\eta(C_n)=n$ holds trivially.
For general groups, an elementary construction gives
\begin{equation}\label{intro:lower}
 \sconst(G)\ge\eta(G)+n-1.
\end{equation}
Gao \cite[Conjecture~2.3]{Gao03} conjectured that equality always holds
in~\eqref{intro:lower}.

Several substantial classes of groups satisfy this equality.
For $1\le m\mid n$, the exact formulas are
\[
 \eta(C_m\oplus C_n)=2m+n-2,\qquad
 \sconst(C_m\oplus C_n)=2m+2n-3;
\]
see~\cite[Theorem~6.3]{GG06}. The determination of the second
constant rests on Reiher's proof of Kemnitz' conjecture~\cite{Reiher07}.
Gao proved the equality for every group of exponent $2$, $3$ or
$4$~\cite[Theorem~2.5]{Gao03}.
There are also results beyond rank two: Girard and
Schmid~\cite[Theorems~3.1 and~3.2]{GS19}, for example, established
it for $C_2\oplus C_{2m}\oplus C_{2m}$ for every $m\ge1$.

In this paper, we prove that this conjectured equality in~\eqref{intro:lower} fails for the group $C_2^7\oplus C_6$ as given in the following theorem.

\begin{theorem}\label{thm:main}
For $G_\ast=C_2^7\oplus C_6$, we have $\sconst(G_\ast)>\eta(G_\ast)+\exp(G_\ast)-1$.
\end{theorem}

\section{Notation and preliminaries}\label{sec:prelim}

We follow the conventions of~\cite[Section~2]{GG06}.
The free abelian monoid over $G$ is denoted by $\mathcal F(G)$.
Its elements are sequences, written multiplicatively:
\[
 S=g_1\cdots g_\ell=\prod_{g\in G}g^{v_g(S)}.
\]
Here $|S|=\ell$ is the length, $v_g(S)$ is the multiplicity of $g$,
$\supp(S)=\{g:v_g(S)>0\}$ is the support, and
$\sigma(S)=\sum_{i=1}^{\ell}g_i$ is the sum.
The empty sequence is denoted by $1$.
We write $T\mid S$ when $T$ is a subsequence of $S$; its complementary
subsequence is $T^{-1}S$.
A sequence is \emph{squarefree} if all multiplicities are at most one.
We identify a squarefree sequence with its support whenever convenient.

For $k\ge1$, put
\[
 \Sigma_k(S)=\{\sigma(T):T\mid S,\ |T|=k\},\quad
 \Sigma_{\le k}(S)=\bigcup_{j=1}^k\Sigma_j(S),\quad
 \Sigma(S)=\bigcup_{j\ge1}\Sigma_j(S).
\]
Thus the empty subsequence is excluded from $\Sigma(S)$ and
$\Sigma_{\le k}(S)$. When it is to be included, we write
\[
 \Sigma_0(S)=\{0\},\qquad
 \Sigma_{\le k}^{0}(S)=\{0\}\cup\Sigma_{\le k}(S).
\]
For a subset $X\subset G$, these definitions refer to the squarefree
sequence $\prod_{x\in X}x$.
In particular, $\Sigma_2(X)$ consists of sums of \emph{distinct}
elements, whereas $X+X$ allows the two summands to coincide.

A sequence $S$ is \emph{zero-sumfree} if $0\notin\Sigma(S)$.

The Davenport constant $\Dconst(G)$ is the least integer $\ell$
such that every sequence over $G$ of length at least $\ell$
has a nonempty zero-sum subsequence.

A homomorphism $\varphi:G\to G'$ is applied to sequences term by term.
In particular, $\varphi(S)=\varphi(g_1)\cdots\varphi(g_\ell)$ and
$\sigma(\varphi(S))=\varphi(\sigma(S))$.
Translation also acts term by term: $h+S=(h+g_1)\cdots(h+g_\ell)$.

We identify $C_2^r$ with the vector space $\F_2^r$.
Let $X=\{x_1,\ldots,x_N\}\subset\F_2^r$ be nonempty.
Its even relation code is
\[
 K_X=\left\{\varepsilon\in\F_2^N:
       \sum_i\varepsilon_i=0,\quad\sum_i\varepsilon_ix_i=0\right\}.
\]
The Hamming weight $\wt(\varepsilon)$ is the number of nonzero coordinates.
Thus a word of weight $2j$ in $K_X$ corresponds to a zero-sum
subsequence of $X$ of length $2j$.
The matrix with columns $(1,x_i)$ has rank
$1+\dim\operatorname{aff}(X)$, so
\[
 \dim K_X=N-1-\dim\operatorname{aff}(X).
\]
Let $T_X\le K_X$ be generated by the words of weight four, called
\emph{tetrads}. Their supports connect the coordinates into components.
A coordinate lying in no tetrad is called \emph{unused}; otherwise it
is called \emph{active}. The notation $0_u$ means $u$ unused coordinates.
These coding terms will be used only to describe the possible
four-term zero-sum relations of a squarefree sequence.

We repeatedly use the following elementary observation. For two binary
words $x,y$,
\[
 \wt(x)+\wt(y)+\wt(x+y)
 =2|\supp(x)\cup\supp(y)|.
\]
Consequently, if three nonzero words $x,y,x+y$ each have weight at
least $m$ and are supported on $N$ coordinates, then $3m\le2N$.
This will rule out several small relation spaces.

\section{Proof of Theorem \ref{thm:main}}\label{sec:reduction}

In this section, we shall prove Theorem~\ref{thm:main}. Throughout this section, identify
\[
 G_\ast=C_2^7\oplus C_6\cong V\oplus C_3,
 \qquad V=C_2^8=\F_2^8,
\]
and let $\pi:G_\ast\to V$ and $\chi:G_\ast\to C_3$ be the projections.
We write the elements of $C_3$ as $0,1,-1$.
For a sequence $S$, let $S_i$ be its subsequence of all terms with
$C_3$-coordinate $i$. Thus $S=S_0S_1S_{-1}$.
Whenever $\pi(S)$ is squarefree, we put
\[
 A=\supp(\pi(S_0)),\quad B=\supp(\pi(S_1)),\quad
 C=\supp(\pi(S_{-1})),\qquad (a,b,c)=(|A|,|B|,|C|).
\]
The sets $A,B,C$ are then pairwise disjoint. We use subsets of $V$
only after making this identification with squarefree sequences.

The rest of this section contains three subsections.

\subsection{The lower bound and the reduction to squarefree projections}
\begin{lemma}\label{reduction:lower}
There is a sequence over $G_\ast$ of length $32$ without a zero-sum
subsequence of length $6$.
\end{lemma}
\begin{proof}
Let $e_1,\ldots,e_6$ be a basis of $\F_2^6$, and set
\[
 C_0=\{0,e_1,\ldots,e_6,e_1+\cdots+e_6\},\qquad
 B_0=C_0\times\F_2\subset\F_2^7.
\]
The only nonempty subset of $C_0$ having even cardinality and sum
zero is $C_0$ itself. Indeed, a relation involving
$e_1+\cdots+e_6$ must involve every $e_i$, and the point $0$ is then
needed to make its cardinality even.
Suppose that six distinct points of $B_0$ had sum zero. The points
of $C_0$ whose fibers contribute an odd number of terms would form
an even zero-sum subset of size at most six. That subset is empty,
so the six points must be three complete fibers of size two.
Their last coordinates sum to $1$, a contradiction.

Now use $G_\ast=\F_2^7\oplus C_6$, and take
\[
 S=\prod_{b\in B_0}(b,0)(b,1).
\]
If $T\mid S$ has length six and sum zero, its $C_6$-coordinate forces
all six terms to have coordinate $0$ or all six to have coordinate $1$.
Their first coordinates would be six distinct elements of $B_0$
summing to zero, which is impossible. Hence $\sconst(G_\ast)\ge33$.
This is the binary doubling construction of~\cite[Theorem~5.4]{Sid20}.
\end{proof}

For the upper bound, suppose that $S\in\mathcal F(G_\ast)$ has no
short zero-sum subsequence. Three elementary restrictions will be useful.
First, $|S_0|\le9$. Otherwise ten of its projected terms have two
independent binary relations $x,y$. Every nonzero relation has
weight at least seven, contradicting $21\le20$ for $x,y,x+y$.
Second, $v_0(\pi(S))\le2$, since $\Dconst(C_3)=3$.
Third, $S$ has at most two pairwise disjoint subsequences $P$ of
length two with $\sigma(\pi(P))=0$. Such a $P$ has
$\chi(\sigma(P))\ne0$. Two disjoint such subsequences cannot have
opposite sums in $C_3$, and three with the same sum would together
form a zero-sum subsequence of length six.

We reduce the remaining argument to the following two assertions.

\medskip\noindent
\textbf{(R) Residual subsequences.}
Let $R\in\mathcal F(G_\ast)$ satisfy
\begin{enumerate}
\item $\pi(R)$ is squarefree, $0\notin\supp(\pi(R))$, and $|R_0|\le9$;
\item if $T\mid R$, $\sigma(\pi(T))=0$ and $|T|\in\{3,4\}$,
      then $\chi(\sigma(T))=1$;
\item if $T\mid R$, $\sigma(\pi(T))=0$ and $|T|\in\{5,6\}$,
      then $\chi(\sigma(T))\ne0$.
\end{enumerate}
Then $|R|\le23$, and equality implies $R_{-1}=1$.

\medskip\noindent
\textbf{(P) Squarefree projections.}
If $S\in\mathcal F(G_\ast)$ has length $27$, $\pi(S)$ is squarefree,
$0\notin\supp(\pi(S))$ and $|S_0|\le9$, then $S$ has a short
zero-sum subsequence.
\medskip

\begin{lemma}\label{reduction:RP}
Assertions \textup{(R)} and \textup{(P)} imply $\eta(G_\ast)\le27$.
\end{lemma}
\begin{proof}
Suppose that $|S|=27$ and $S$ has no short zero-sum subsequence.
If $\pi(S)$ is squarefree and does not contain $0$, assertion (P) applies.
Otherwise choose $E\mid S$ as follows: take two terms with equal
projection if such terms exist, and otherwise take one term
projecting to zero. Thus $|E|\in\{1,2\}$ and $\sigma(\pi(E))=0$.
Its $C_3$-sum is nonzero. Applying the automorphism
$(v,i)\mapsto(v,-i)$ if necessary, assume $\chi(\sigma(E))=1$.

Every subsequence $R\mid E^{-1}S$ whose projection is squarefree
and avoids $0$ satisfies the hypotheses of (R). To see the only
additional condition, let $T\mid R$ have projected sum zero and
length three or four. Its $C_3$-sum cannot be $0$, since $T$ is short,
and cannot be $-1$, since $ET$ would be short and zero-sum.
It must therefore be $1$. The conditions for lengths five and six
follow directly from the assumption on $S$.

Suppose first that $S$ contains no two disjoint pairs of terms
with equal projection. Delete the chosen subsequence $E$.
After also removing any
remaining term projecting to zero, at least $24$ terms remain,
with distinct nonzero projections. They give a contradiction to (R).
Here a third occurrence of a projected value causes no difficulty:
removing two occurrences leaves only one; four occurrences or two
different repeated values would give two disjoint pairs.

It remains to consider two disjoint pairs $P,Q\mid S$, each with
projected sum zero. Their $C_3$-sums agree, and we may normalize both
to $1$. Put $R=(PQ)^{-1}S$, so $|R|=23$.
There is no repeated projection in $R$, since it would give a third
disjoint pair. Nor can a term $g\mid R$ project to zero:
according as $\chi(g)=0,-1,1$, one of $g$, $Pg$, $PQg$ would be short
and zero-sum. Assertion (R) therefore applies, and all terms of $R$
have $C_3$-coordinate $0$ or $1$.

Let $x$ be the common projection of the terms of $Q$. If $x\ne0$
and $x\notin\supp(\pi(R))$, adjoining either term of $Q$ to $R$
gives a sequence of length $24$ satisfying (R), using $P$ as the
deleted subsequence $E$. If $x\in\supp(\pi(R))$, the two terms of
$Q$ and the term of $R$ with projection $x$ give three possible
pairs, each disjoint from $P$. Their $C_3$-sums must all be $1$.
Writing their coordinates as $i,j,k$, the equations
$i+j=i+k=j+k=1$ force $i=j=k=-1$, contrary to the conclusion for $R$.
Thus $x=0$. The same argument applies to $P$, yielding four terms
projecting to zero, contrary to $v_0(\pi(S))\le2$.
\end{proof}

We now prove (R) and (P). The next subsection establishes the common
counting estimates and proves (R); the last subsection proves (P)
by counting in three-dimensional quotients.

\subsection{Binary structure and residual subsequences}\label{rr:section}

We first count three-term sums of a squarefree sequence in $V$ with
no zero-sum subsequence of length six. The four-term relations determine
this count. After establishing the counting and overlap estimates, we
apply them to the projected coordinate classes and prove assertion (R).

\begin{lemma}\label{rr:blocks}
Let $X\subset V$ satisfy $0\notin\Sigma_6(X)$. The code $T_X$ is a
direct sum of unused coordinates and blocks of the following types:
\begin{enumerate}
\item $d_{2m}$, for $m\ge2$: partition the $2m$ coordinates into
      $m$ pairs, and take the unions of an even number of these pairs;
\item $e_7$ and $e_8$: the even relation codes of a punctured affine
      three-space and a full affine three-space, respectively.
\end{enumerate}
The pairs in a $d_{2m}$ block all have the same nonzero sum.
\end{lemma}
\begin{proof}
Two tetrads cannot meet in one point: their symmetric difference
would be a zero-sum subsequence of length six. They cannot meet in
three points either, since their symmetric difference would have
length two and $X$ is squarefree. Thus any two tetrads meet evenly.
The decomposition theorem of Pless--Sloane
\cite[Theorem~6.5]{PlessSloane} gives precisely the stated blocks.

In $d_{2m}$, the union of any two coordinate pairs is a tetrad.
Consequently the two pairs have the same sum; denote this sum by $d$.
The corresponding points are therefore pairs $\{x,x+d\}$.
An $e_s$ block, $s\in\{7,8\}$, supplies $s-4$ independent even
relations, so its points have affine dimension at most three.
There are more than four distinct points, so the affine dimension
is exactly three. This identifies the points with the stated
punctured or full affine three-space.
\end{proof}

We record two geometric consequences that will be used later.
First, suppose that $X$ consists of two $e_8$ blocks, written
$x+U_0$ and $y+W_0$, where $U_0,W_0$ have dimension three.
Their direction spaces meet trivially: a common nonzero direction
would give a tetrad meeting both blocks. Also $x+y\notin U_0+W_0$,
because the two affine spaces are disjoint. Hence
\begin{equation}\label{rr:eight-geometry}
 L=x+U_0+W_0+\langle x+y\rangle
 \quad\text{has affine dimension seven},\qquad
 \Sigma_3(X)=L,\quad \Sigma_4(X)=L+L.
\end{equation}
For the three-term sums, separate the selections according to
which block contributes an odd number of points. The selections
using one point of the first block and two of the second, together
with those using three of the first, fill $x+U_0+W_0$; reversing
the blocks gives $y+U_0+W_0$. For four-term sums,
selections of types $4+0$, $2+2$ and $0+4$ give $U_0+W_0$;
selections of types $3+1$ and $1+3$ give $x+y+U_0+W_0$.
These are exactly the sets in \eqref{rr:eight-geometry}.

Second, let $D$ consist of $m\ge7$ pairs with common sum $d$, and
write $q_1,\ldots,q_m$ for their distinct images in $V/\langle d\rangle$.
These indices have no even relation of length two, four or six.
A four-term relation, together with a spare complete pair, would
give six points summing to zero: choose the four lifts to sum to
$d$, and add the spare pair. A six-term relation gives the same
contradiction by choosing lifts with sum zero. In particular, any
seven indices are affinely independent, and the corresponding
fourteen points span an affine seven-space.

The following description of the sumsets is useful. In a selection
of points of $D$, record the indices used an odd number of times.
Changing the chosen point at any such index changes the sum by $d$.
It follows that $\Sigma_3(D)$ consists of both lifts of every
one-index and three-index sum. Similarly, $\Sigma_4(D)$ consists
of both lifts of the two-index and four-index sums, together with
$0$, obtained from two complete pairs. The available spare pairs
ensure that all these representations use distinct points.

\begin{lemma}\label{rr:cosets}
For a direct sum $T$ of the blocks above, let $b(T)$ and $r(T)$ be
the numbers of cosets of $T$ represented by words of weight two and
three, respectively. The component counts are
\begin{equation}\label{rr:component-counts}
\begin{array}{c|c|c|c}
T&\text{length}&b(T)&r(T)\\\hline
d_{2m}&2m&1+2\binom m2&2m+2\binom m3\\
e_7&7&7&8\\
e_8&8&7&8\\
0_1&1&0&0
\end{array}
\end{equation}
For disjoint components of lengths $n_1,n_2$,
\begin{align}
b(T_1\oplus T_2)&=b(T_1)+b(T_2)+n_1n_2,\label{rr:b-sum}\\
r(T_1\oplus T_2)&=r(T_1)+r(T_2)+b(T_1)n_2+n_1b(T_2).
\label{rr:r-sum}
\end{align}
If $0\notin\Sigma_6(X)$, then
\begin{equation}\label{rr:triple-count}
 r(T_X)=|\Sigma_3(X)|\le2^{\dim\operatorname{aff}(X)}.
\end{equation}
\end{lemma}
\begin{proof}
For $d_{2m}$, a weight-two word either selects one complete pair,
giving one coset, or selects two different pair indices, giving two
cosets for each pair of indices. A weight-three word has either
one or three indices used oddly. These give $2m$ and
$2\binom m3$ cosets, respectively.

For $e_7,e_8$, every nonzero direction in the affine three-space
is a pair sum, giving seven possibilities. Every point of that
affine space is also a three-term sum. Indeed, for a target $t$,
choose a point $x\ne t$ of the block. There are at least three pairs
with sum $t+x$; at most one contains $x$, so one of the others
completes the required triple. This gives eight possibilities.

For a direct sum, distribute the selected coordinates between
the two components. Weight two has distributions $(2,0),(1,1),(0,2)$,
and weight three has distributions $(3,0),(2,1),(1,2),(0,3)$.
They give the terms in \eqref{rr:b-sum} and \eqref{rr:r-sum}.
The relevant cosets are distinct, since the component codes are
even and have no words of weight two.

Finally, two distinct three-term selections of $X$ with equal
sum have symmetric difference of weight two, four or six.
Only weight four is possible, and this difference is a tetrad.
Thus their incidence vectors lie in the same coset of $T_X$,
and conversely a common coset gives equal sums. All odd sums
belong to the affine span of $X$, proving \eqref{rr:triple-count}.
\end{proof}

For a block code $T$ of length $n$, let $d(T)$ be the number of
blocks $d_{2m}$ and let $u(T)$ be the number of unused coordinates.
Define
\[
 F_n(T)=r(T)-n\,d(T)-(n-2)\left\lfloor\frac{u(T)}2\right\rfloor,
 \qquad g_0(n)=2\binom{\lceil n/2\rceil}{3},
\]
and put $g(n)=g_0(n)$ except that $g(9)=15$.
As usual, a binomial coefficient is zero when its upper argument
is smaller than its lower argument. The subtractions in $F_n$
will account for overlaps between three-term sums and translates.

\begin{lemma}\label{rr:counting}
For $n\le19$, one has $F_n(T)\ge g(n)$. Moreover,
\begin{enumerate}
\item $F_{14}(T)\le99$ forces $T=d_{14}$;
\item $F_{16}(T)\le153$ forces $T=d_{16}$ or $e_8\oplus e_8$,
      with values $112$ and $128$, respectively;
\item $F_{18}(T)\le223$ forces $T=d_{18}$.
\end{enumerate}
Every squarefree sequence over $V$ of length nineteen has a
zero-sum subsequence of length six.
\end{lemma}
\begin{proof}
We reduce the calculation to a single paired block and at most
two exceptional blocks $e_7,e_8$.
Replace each pair of unused coordinates formally by a block
$d_2$ with counts $(n,b,r)=(2,1,2)$. This is only a counting device.
It increases $r(T)$ by two and increases the subtracted quantity
in $F_n(T)$ by two, so $F_n$ is unchanged.

Next merge paired blocks of half-lengths $x,y$ into $d_{2(x+y)}$.
Using \eqref{rr:b-sum} and \eqref{rr:r-sum}, the decrease in $F_n$ is
\[
 xy(2n-2-3x-3y)\ge0.
\]
Indeed, $n\ge2(x+y)$ and $x+y\ge2$. Thus repeated merging gives
a lower bound for the original $F_n$.

Call the sum of the $e_7,e_8$ blocks the exceptional core, and let
its length, pair count and triple count be $(s,p,r)$. Write $m=n-s$.
After merging, the remaining coordinates consist of one paired
block and, when $m$ is odd, one unused coordinate. Applying the
direct-sum formulas gives
\begin{equation}\label{rr:core-bound}
 F_n(T)\ge Q_n(s,p,r):=
 r+pm+2s\binom{\lceil m/2\rceil}{2}+g_0(m).
\end{equation}
For example, when $m=2\ell$, the paired block has counts
$b=1+2\binom\ell2$, $r=m+2\binom\ell3$; its subtraction of $n=s+m$
cancels the terms $s+m$ in the direct-sum formula. For odd $m$,
the same calculation uses
$\binom\ell3+\binom\ell2=\binom{\ell+1}3$.
The cases $m=0,1$ follow directly as well.

An empty core gives $Q_n=g_0(n)$. For one exceptional block,
the excess $Q_n-g_0(n)$ is as follows:
\[
\begin{array}{c|cc}
\text{core}&m=2\ell&m=2\ell+1\\\hline
e_7&\ell(3\ell-1)&4\ell^2+12\ell+7\\
e_8&2\ell(2\ell-1)&4\ell^2+6\ell-5
\end{array}
\]
Every entry is nonnegative except the last one at $\ell=0$.
That exception is $e_8\oplus0_1$, whose value is $F_9=15$,
explaining the definition of $g(9)$.

There can be at most two exceptional blocks, since $n\le19$.
Their core counts, obtained from \eqref{rr:b-sum}--\eqref{rr:r-sum}, are
\[
\begin{array}{c|ccc}
\text{core}&s&p&r\\\hline
e_7\oplus e_7&14&63&114\\
e_7\oplus e_8&15&70&121\\
e_8\oplus e_8&16&78&128
\end{array}
\]
In particular $Q_n\ge114+63m$. For $m=0$, we have $n\le16$
and $g_0(n)\le112$; for $m=1$, we have $n\le17$ and $g_0(n)\le168$;
for $m\ge2$, we have $g_0(n)\le g_0(19)=240\le114+63m$.
This proves $F_n\ge g(n)$; lengths at most three can also be
checked immediately from the definitions.

We now retain the information needed for the three small profiles.
For a core-free code of even length $n=2\ell$ other than $d_n$,
the last merger has half-lengths $x+y=\ell$ and decreases $F_n$ by
$xy(\ell-2)\ge(\ell-1)(\ell-2)$. Since $F_n(d_n)=g_0(n)$,
the resulting lower bounds at $n=14,16,18$ are $100,154,224$.
For nonempty cores, \eqref{rr:core-bound} gives the following table:
\[
\begin{array}{c|rrr}
\text{core}&Q_{14}&Q_{16}&Q_{18}\\\hline
e_7&149&231&335\\
e_8&100&168&258\\
e_7\oplus e_7&114&240&394\\
e_7\oplus e_8&--&191&361\\
e_8\oplus e_8&--&128&284
\end{array}
\]
A dash means the core is too long. Comparing these entries with
$99,153,223$ proves the asserted profiles.

Finally, a block code of length nineteen has a paired block or at
least three unused coordinates: without either, the two possible
exceptional blocks and at most two unused coordinates use at most
eighteen coordinates. Therefore
\[
 r(T)\ge F_{19}(T)+17\ge g(19)+17=257.
\]
For a squarefree sequence in $V$ with $0\notin\Sigma_6(X)$ this
contradicts $r(T_X)=|\Sigma_3(X)|\le256$.
\end{proof}

\begin{lemma}\label{rr:overlap}
Suppose $B\subset V$ satisfies $0\notin\Sigma_6(B)$, and let
$E\subset V$ satisfy
\[
 E\cap\Sigma_4(B)=\varnothing,
 \qquad \Sigma_2(E)\cap\Sigma_2(B)=\varnothing.
\]
Then
\[
 |\Sigma_3(B)\setminus(B+E)|\ge F_{|B|}(T_B).
\]
\end{lemma}
\begin{proof}
Put $n=|B|$. If a sum of three distinct points of $B$ equals
$b+e$, with $b\in B$ and $e\in E$, then that triple contains $b$.
Otherwise the four distinct $B$-points would represent $e$,
contrary to $E\cap\Sigma_4(B)=\varnothing$.
The other two points therefore form a pair representing $e$.

Such a pair meets every tetrad evenly. If it met a tetrad in
exactly one point, their symmetric difference would be a
four-point representation of $e$. It follows from the block
description that the pair is either entirely unused or lies in
one paired block. Indeed, an active pair crossing two blocks
meets some tetrad in one point. In $e_7$ or $e_8$, an affine plane
can be chosen through one endpoint and avoiding the other (and
the missing point in the $e_7$ case), with the same consequence.
Within $d_{2m}$, $m\ge3$, a third coordinate pair shows that a pair
meeting every tetrad evenly must be one of the distinguished
pairs, so its sum is the common direction. For $d_4$, there are
three possible pair directions, but two of them cannot both
belong to $E$: their sum belongs to $\Sigma_2(B)$.

Pairs of unused coordinates have unique representations, since
two different pairs with the same sum would produce a tetrad.
Two selected unused pairs cannot share an endpoint either:
the sum of their two different represented elements of $E$
would then belong to $\Sigma_2(B)$.

Thus at most $d(T_B)$ block directions can contribute to the
overlap, each through at most $n$ choices of $b$. At most
$\lfloor u(T_B)/2\rfloor$ unused pairs contribute, each through
at most $n-2$ choices, since $b$ must avoid that pair. Subtracting
these upper bounds from $r(T_B)=|\Sigma_3(B)|$ proves the claim.
\end{proof}

\begin{proposition}\label{rr:mixed}
Suppose $S\in\mathcal F(G_\ast)$ has no short zero-sum subsequence,
$\pi(S)$ is squarefree, and $0\notin\supp(\pi(S))$.
Use the classes $A,B,C$ and lengths $a,b,c$ defined above, and put
$H(a)=1+a+\binom a2+\binom a3$. Then
\begin{equation}\label{rr:mixed-bound}
 H(a)+bc+F_b(T_B)+F_c(T_C)\le256,
 \qquad b+c\le24.
\end{equation}
If $|S|=27$, then $3\le a\le8$ and $b,c\le16$.
\end{proposition}
\begin{proof}
All sums of at most three distinct points of $A$, including the
empty sum, are different: a collision gives a nonempty zero-sum
subsequence of $S_0$ of length at most six. Thus there are $H(a)$
such sums. Also $B+C$ has $bc$ elements, because two different
representations give a zero-sum subsequence with two terms of
each nonzero $C_3$-coordinate. The four sets
\[
 \Sigma_{\le 3}^{0}(A),\quad B+C,\quad
 \Sigma_3(B)\setminus(B+C),\quad
 \Sigma_3(C)\setminus(B+C)
\]
are pairwise disjoint: an intersection would again give a short
zero-sum subsequence. For example, three terms from $B$ have
$C_3$-sum zero, as do three from $C$.
Lemma~\ref{rr:overlap} applies with $(B,E)=(B,C)$, and also with
$B,C$ interchanged. Its hypotheses follow by excluding projected
relations with four terms from one class and one from the other,
or two from each. This proves the first inequality.

By Lemma~\ref{rr:counting}, $b,c\le18$. If $b+c\ge25$, take a
subsequence with $b+c=25$ and write its even class size as $2k$.
Then $4\le k\le9$. On the even class $g=g_0$, and on the odd
class $g\ge g_0-5$, so expansion gives
\[
 bc+g(b)+g(c)\ge259+(k-7)(7k-44)\ge259.
\]
For the six possible integers $k$, the last product is nonnegative.
This contradicts the first inequality and proves $b+c\le24$.

Now suppose $|S|=27$. The preliminary bound $a\le9$ and $b+c\le24$
give $3\le a\le9$. For these integers,
\[
 H(a)-16a+54=\frac{(a-5)(a-6)(a+11)}6\ge0.
\]
If $b\in\{17,18\}$, then $bc\ge18(9-a)$ and $F_b(T_B)\ge168$.
Hence \eqref{rr:mixed-bound} implies
\[
 256\ge H(a)+bc+168\ge276-2a\ge258,
\]
a contradiction. The same argument applies to $c$.

It remains to exclude $a=9$. The binary relation space of $A$ has
dimension one: it is nonzero because $A\subset V$, and two
independent relations would give $21\le18$, since every nonzero
relation has weight at least seven. Its unique nonzero word has
weight seven, eight or nine. Among the $256$ selections of at
most four points of $A$, a repeated sum comes from complementary
parts of that word. There are at most $35$ pairs of equal sums
($\binom73=35$ or $\binom84/2=35$). Thus
$|\Sigma_{\le 4}^{0}(A)|\ge221$. This set is disjoint from $B+C$, so
$bc\le35$. Since $b+c=18$ and $b,c\le16$, only
$\{b,c\}=\{16,2\}$ remains. But then
\[
 256\ge H(9)+32+g(16)=130+32+112=274,
\]
the final contradiction.
\end{proof}

\begin{lemma}\label{rr:two-coordinates}
Let $R$ satisfy the hypotheses of \textup{(R)}. If
$\supp(\chi(R))\subseteq\{0,1\}$, then $|R|\le23$.
\end{lemma}
\begin{proof}
The hypotheses are inherited by subsequences, so it suffices to
exclude $|R|=24$. Let $A,B$ be the projected classes with
$C_3$-coordinates $0,1$, of lengths $a,b$. The preliminary bound
$a\le9$ and Lemma~\ref{rr:counting}, which gives $b\le18$, imply
$6\le a\le9$. Put
\[
 A^0=\{0\}\cup A,\qquad Y=B+A^0,\qquad M=|Y|.
\]
The translates $b_0+A^0$, for $b_0\in B$, are pairwise disjoint.
Indeed, a collision between two translates gives two distinct
points of $B$ together with one or two points of $A$, having
projected sum zero and $C_3$-sum $2$, contrary to (R).
Also $0\notin Y$, because $A^0\cap B=\varnothing$. Consequently
\begin{equation}\label{rr:Y-size}
 M=b(a+1).
\end{equation}

Let $k$ be the number of three-element subsets of $A^0$ whose sums
belong to $B$. Their sums are all different: two representations
would give a nonempty zero-sum subsequence of $A$ of length at
most six, forbidden by the hypotheses on $R$.
Every represented point of $B$ is unused in $T_B$. Otherwise a
tetrad through that point, combined with its representation by
$A^0$, gives three points of $B$ and two or three points of $A$
with total projected sum zero and $C_3$-sum zero.

There cannot be eight unused coordinates of $T_B$. Choose eight
of them and any seven other coordinates, which exist because
$b\ge15$. Three-term selections supported on these fifteen
coordinates give at least
\begin{equation}\label{rr:unused-eight}
 \binom83+7\binom82+8\cdot6=56+196+48=300
\end{equation}
distinct cosets of $T_B$. The three terms count selections with
three, two and one unused coordinates. For the last term, the
other seven coordinates have at least six distinct pair sums:
fix one of them and pair it with each of the other six.
Different subsets of unused coordinates cannot be identified
modulo $T_B$. Thus \eqref{rr:unused-eight} contradicts
$r(T_B)=|\Sigma_3(B)|\le256$. We have proved
\begin{equation}\label{rr:k-bound}
 k\le7.
\end{equation}

We now count ordered zero-sum triples in $Y$; write their number
as $\tau$. Their terms are distinct because $0\notin Y$.
Each $y\in Y$ has a unique expression $b_0+a_0$ with
$b_0\in B$, $a_0\in A^0$. In a zero-sum triple, the three
$B$-coordinates cannot all be distinct: together with the
uncancelled $A$-coordinates they would give a forbidden subsequence
of $R$ of length between three and six and $C_3$-sum zero.
The three $A^0$-coordinates are distinct as well. If two were
equal, cancelling them and the repeated $B$-coordinate would
force a point of $B$ to belong to $A^0$.

To bound $\tau$ from above, choose the pair of positions with
equal $B$-coordinates, their common value, and the ordered triple
from $A^0$. There are at most $3$, $b$ and $6k$ choices, respectively;
the remaining $B$-coordinate is then determined by the zero-sum
condition. This may overcount, which is harmless. Hence
\[
 \tau\le18bk.
\]
For the lower bound, fix the first term $y\in Y$. The possible
second terms form $Y\cap(y+Y)$, and each determines the third.
Since both sets have size $M$ in $V$,
$|Y\cap(y+Y)|\ge2M-256$. Thus
\begin{equation}\label{rr:triangles}
 M(2M-256)\le\tau\le18bk.
\end{equation}
For $a\ge7$, \eqref{rr:Y-size} gives $M\ge136$, while $b\le17$.
Together with \eqref{rr:k-bound}, this gives
$2176\le\tau\le18\cdot17\cdot7=2142$, a contradiction.

Therefore $a=6$ and $b=18$. As in the proof of
Proposition~\ref{rr:mixed}, the sets $\Sigma_{\le 3}^{0}(A)$ and
$\Sigma_3(B)$ are disjoint. Hence $r(T_B)\le256-H(6)=214$.
In particular $F_{18}(T_B)\le214$, so Lemma~\ref{rr:counting}
forces $T_B=d_{18}$. There are no unused coordinates, hence $k=0$.
By \eqref{rr:triangles}, $Y$ has size $126$ and is sum-free,
meaning that $Y\cap(Y+Y)=\varnothing$.

We use the Davydov--Tombak classification in the form recorded
in~\cite[Fact~2]{GL10}. A maximal sum-free subset of $\F_2^r$
of size greater than $9\cdot2^{r-5}$ is either a nonzero
hyperplane coset or has size $5\cdot2^{r-4}$.
Extend $Y$ to a maximal sum-free set. Since $126>80=5\cdot2^4$,
only the hyperplane-coset alternative is possible. As $B\subseteq Y$,
the nine parallel pairs of $B$ lie in that coset.
Their common direction belongs to its direction hyperplane;
after quotienting by this direction, their nine indices have
affine dimension at most six. Their even relation space therefore
has dimension at least two. Every nonzero even relation has
weight at least eight, by the parallel-pair description above.
Two independent relations and their sum would then give
$24\le2\cdot9=18$, the final contradiction.
\end{proof}

\begin{proposition}\label{rr:residual}
Assertion \textup{(R)} holds.
\end{proposition}
\begin{proof}
First let $|R|=23$ and assume the hypotheses of (R). Adjoin the
term $(0,0)$ and translate every term by $(0,-1)$; call the
resulting sequence $S$. Its projection is squarefree. For every
$Q\mid S$ with $\sigma(\pi(Q))=0$, we have
\begin{equation}\label{rr:translated-conditions}
 \chi(\sigma(Q))=0\quad\text{if }|Q|=4,
 \qquad
 \chi(\sigma(Q))\ne0\quad\text{if }|Q|=6.
\end{equation}
For a four-term subsequence not using the adjoined term, the
translation subtracts $4\equiv1$ from its original $C_3$-sum $1$.
If it uses the adjoined term, the other three terms originally
had $C_3$-sum $1$, so the same conclusion holds. For length six,
the translation does not change the sum; the original five- or
six-term condition applies according as the adjoined term occurs.

Let $U,B,W\subset V$ be the projected classes with $C_3$-coordinates
$-1,0,1$, and put $(u,b,w)=(|U|,|B|,|W|)$. They are pairwise
disjoint and $u+b+w=24$. Suppose, towards a contradiction, that
$u,w>0$. The conditions \eqref{rr:translated-conditions} are invariant
under the automorphism $(v,i)\mapsto(v,-i)$, so we may assume $u\ge w$.

\emph{Step 1: disjoint sumsets.}
Neither $U$ nor $W$ has a zero-sum subsequence of length four or
six: four equal nonzero $C_3$-coordinates have nonzero sum, whereas
six have sum zero, in both cases contradicting
\eqref{rr:translated-conditions}. Consequently their two-term sums
have unique unordered representations, and all their one- and
three-term sums are distinct. Moreover $u,w\le10$. Eleven points
in $V$ have at least two independent even relations, and the
three nonzero words they generate would each have weight at
least eight, contradicting $24\le22$.

Put
\[
 T=U+W,\qquad F_U=U\cup\Sigma_3(U),\qquad
 F_W=W\cup\Sigma_3(W).
\]
The sets $F_U,F_W$ and all translates $x+T$, $x\in B$, are
pairwise disjoint. Here is how \eqref{rr:translated-conditions}
rules out the possible intersections.
\begin{itemize}
\item An intersection of $F_U$ and $F_W$ gives either a four-term
relation with $C_3$-sum $1$ or $-1$, or a six-term relation with
$C_3$-sum zero. The one-term sets themselves are disjoint.
\item An intersection of $F_U$ with $x+T$ gives, after cancellation,
coordinates $0,1,-1,-1$, or, without cancellation in the
three-term case, $0,1,-1,-1,-1,-1$. The sums are respectively
nonzero for length four and zero for length six. A cancellation
leaving only two terms is impossible because the projected
classes are disjoint. The argument for $F_W$ is symmetric.
\item If $x+T$ and $x'+T$ meet with $x\ne x'$, cancellation leaves
two terms from $B$ and two from one nonzero coordinate class,
or two from each of the three classes. These are again forbidden
relations of length four or six.
\end{itemize}
We also have
\begin{equation}\label{rr:additional-disjointness}
 \Sigma_3(B)\cap(F_U\cup F_W)=\varnothing,
 \qquad T\cap\Sigma_4(B)=\varnothing.
\end{equation}
The first exclusion uses three terms from $B$ together with one
or three from $U$ or $W$; the second uses four from $B$, one
from $U$ and one from $W$. Their $C_3$-sums are precisely those
forbidden in \eqref{rr:translated-conditions}.
Disjointness of the translates also gives
$\Sigma_2(T)\cap\Sigma_2(B)=\varnothing$.
Applying Lemma~\ref{rr:overlap} with $E=T$, and counting all the
disjoint sets just described, yields
\begin{equation}\label{rr:signed-bound}
 u+w+\binom u3+\binom w3+b|T|+F_b(T_B)\le256.
\end{equation}
Notice also that
\begin{equation}\label{rr:pair-avoidance}
 \bigl(\Sigma_2(U)\cup\Sigma_2(W)\bigr)\cap\Sigma_2(B)=\varnothing:
\end{equation}
both pair-sum sets on the left are contained in $\Sigma_2(T)$,
by translating with one fixed point of the opposite nonempty class.

\emph{Step 2: the possible lengths.}
By \eqref{rr:translated-conditions}, $0\notin\Sigma_6(B)$, so
$b\le18$. Together with $u,w\le10$ this gives
\[
 4\le b\le18,\qquad 6\le u+w\le20.
\]
For $t\in T$, let $r_t$ count its representations as $u_0+w_0$.
Two different representations determine a common element of
$\Sigma_2(U)$ and $\Sigma_2(W)$. Conversely, each common pair
sum gives exactly two unordered pairs of representations.
Since the pair sums in each class are unique,
\[
 \sum_t\binom{r_t}{2}=2|\Sigma_2(U)\cap\Sigma_2(W)|
 \le2\binom w2.
\]
Using $\sum_t r_t=uw$, we obtain
\[
 \sum_t r_t^2\le uw+2w(w-1)=w(u+2w-2).
\]
Cauchy's inequality now gives the single bound
\begin{equation}\label{rr:Cauchy}
 |T|\ge\left\lceil\frac{u^2w}{u+2w-2}\right\rceil.
\end{equation}
Inserting it and $F_b(T_B)\ge g(b)$ into \eqref{rr:signed-bound}
shows that
\[
 L(u,w):=u+w+\binom u3+\binom w3+
 (24-u-w)\left\lceil\frac{u^2w}{u+2w-2}\right\rceil
 +g(24-u-w)\le256.
\]
The following table gives every value in the range
$\max\{w,6-w\}\le u\le10$. A dash denotes a pair outside that range;
the only values at most $256$ are in bold.
\[
\begin{array}{c|rrrrrrrr}
w\backslash u&3&4&5&6&7&8&9&10\\\hline
1&--&--&274&297&267&297&290&331\\
2&--&286&321&284&321&318&360&376\\
3&\mathbf{248}&282&275&307&312&346&365&416\\
4&--&\mathbf{240}&285&272&328&340&383&408\\
5&--&--&\mathbf{240}&280&301&339&368&412\\
6&--&--&--&260&306&330&368&404\\
7&--&--&--&--&294&328&359&404\\
8&--&--&--&--&--&328&368&400\\
9&--&--&--&--&--&--&368&400\\
10&--&--&--&--&--&--&--&404
\end{array}
\]
Thus only $u=w=3,4,5$ remain. Returning to the stronger inequality
\eqref{rr:signed-bound} gives
\[
\begin{array}{c|c|c|c|l}
u=w&b&|T|\ge&F_b(T_B)\le&T_B\\\hline
3&18&4&176&d_{18}\\
4&16&7&128&d_{16}\text{ or }e_8\oplus e_8\\
5&14&10&86&d_{14}
\end{array}
\]
The last column follows from Lemma~\ref{rr:counting}.

\emph{Step 3: the two hyperplane configurations.}
If $T_B=e_8\oplus e_8$, take $D=B$ and use
\eqref{rr:eight-geometry}. If $T_B=d_{18}$, its nine pair indices
lie in the seven-dimensional quotient by their common direction.
They therefore have a nonzero even relation. Its weight must be
eight, since weights two, four and six were excluded above.
The eight indices in that relation give a sixteen-point set $D$.

For completeness, these eight indices have just the one
nonzero even relation, namely their total sum. Their affine span
has dimension six. Their one-index and three-index sums are all
distinct and number $8+\binom83=64$, so they fill an affine
six-space. Their two-index sums and four-index sums, together
with zero, fill its direction space: there are
$\binom82+\binom84/2+1=64$ possibilities, with complementary
four-subsets giving the same sum. Lifting back to $V$ gives
\[
 \Sigma_3(D)=L,\qquad
 |(L+L)\setminus\Sigma_4(D)|=1,
\]
where $L$ is an affine hyperplane. The missing element is the
common pair direction; zero itself is the sum of two complete pairs.
For $e_8\oplus e_8$ the same conclusion holds with no missing element.

By \eqref{rr:additional-disjointness}, both $U$ and $W$ avoid $L$.
The complement of an affine hyperplane is its other coset, so
$U+W\subset L+L$. But $T\cap\Sigma_4(D)=\varnothing$, whence
$|T|\le1$, contradicting the table above.

\emph{Step 4: seven parallel pairs.}
The remaining possibilities are $d_{16}$ and $d_{14}$, with
$u=w=4$ or $5$. Choose seven of their parallel pairs, with common
direction $d$, and let their affine span be $L$.
Translate all projected points so that $0\notin L$. This preserves
every condition concerning sums of four or six points. Project
modulo $\langle d\rangle$ and choose coordinates so that the seven
pair indices are the seven unit vectors of $\F_2^7$.
This is possible because the indices are affinely independent
and their affine span does not contain zero.

The parallel-pair description shows that $\Sigma_3(B)$ contains
both lifts of all words of weights one and three, and
$\Sigma_4(B)$ contains both lifts of all words of weights two and
four. Each of $U,W$ projects injectively: two points with the same
projection differ by $d\in\Sigma_2(B)$, contradicting
\eqref{rr:pair-avoidance}. The same condition excludes a weight-two
difference between two projected points of either class, since
both lifts of any two-index sum belong to $\Sigma_2(B)$.

Separate each projected class into its even-weight and odd-weight
words. An odd word of $U$ or $W$ has weight five or seven, since
weights one and three are excluded by
\eqref{rr:additional-disjointness}. Complementing all seven bits
therefore makes the odd words into zero or two-element subsets.
The two-element subsets in one class must be disjoint: an
intersection in one point would give a weight-two difference.
If the zero word occurs, no two-element subset can occur with it.
Thus each class has at most three odd words.

Since each class has at least four points, each also has an even
word. Fix an even word of the opposite class. Its difference
with any even word of the first class has weight zero or six:
weights two and four would place the sum in $\Sigma_4(B)$.
There is at most one word of each of these two types. Indeed,
two different weight-six words in seven coordinates differ in
weight two, which is forbidden within a class. Hence each class
has at most two even words, and therefore has odd words as well.

Finally, complementing two odd words, one from each class, does
not change their difference. That difference has even weight at
most four, so it must be zero, again because $T$ avoids
$\Sigma_4(B)$. All cross-pairs of odd words are therefore equal;
each odd part has at most one point. Each class would then have
at most one odd and two even words, fewer than four. This is the
required contradiction.

We have proved that one of $U,W$ is empty. Undo any negation of
the $C_3$-coordinate. The adjoined term has coordinate $-1$, so
the class with coordinate $1$ must be empty. Translating back,
every term of the original $R$ has coordinate $0$ or $1$; that is,
$R_{-1}=1$.

If a sequence satisfying the hypotheses of (R) had length at
least $24$, restrict to length $24$. Every one of its terms lies
in a subsequence of length $23$ satisfying the same hypotheses.
The conclusion just proved would give only coordinates $0,1$
throughout, contradicting Lemma~\ref{rr:two-coordinates}.
This proves (R).
\end{proof}

\subsection{The case of distinct binary coordinates}\label{rd:section}
We now prove assertion \textup{(P)}. Suppose, to the contrary, that
$S\in\mathcal F(G_\ast)$ has length $27$, that $\pi(S)$ is squarefree
with $0\notin\supp(\pi(S))$, and that $S$ has no short zero-sum subsequence.
Retain the sets
$A=\supp(\pi(S_0))$, $B=\supp(\pi(S_1))$ and $C=\supp(\pi(S_{-1}))$,
of cardinalities $a,b,c$, respectively.
Proposition~\ref{rr:mixed} gives $3\le a\le8$ and $b,c\le16$;
in particular, $b,c\ge3$.

We shall compare their intersections with cosets of three-dimensional
subspaces. A \emph{coordinate subspace} is a subspace
$H=\langle a_1,a_2,a_3\rangle$ generated by three independent elements of $A$.
For $D\subset V$ and $q\in V/H$, write $D_q=D\cap q$, and put
\[
 h_H(D)=\max_{q\in V/H}|D_q|,
 \qquad
 \nu_H(D)=\sum_{q\in V/H}\min\{2,|D_q|\}.
\]
Thus $h_H(D)$ is the maximum multiplicity of the image of the
squarefree sequence $D$ in $V/H$. The quantity $\nu_H(D)$ counts
one for an intersection of size one and two for an intersection of
size at least two. Empty intersections contribute zero.

\begin{lemma}\label{rd:coset-count}
Let $H$ be a three-dimensional subspace and let $D\subset V$ satisfy
$0\notin\Sigma_6(D)$. At most two $H$-cosets contain three or more
elements of $D$. Moreover,
\[
 |D|\le\nu_H(D)+6,
\]
and the stronger bound $|D|\le\nu_H(D)+5$ holds unless $D$ contains
an entire $H$-coset.
\end{lemma}
\begin{proof}
Three elements in one coset give three pair sums, which are the
nonzero elements of a plane in $H$. If three cosets each contained
at least three elements of $D$, let these planes be $P_1,P_2,P_3$.
We can choose nonzero $\delta_i\in P_i$ with
$\delta_1+\delta_2=\delta_3$. Indeed, if $P_1=P_2$, choose
$\delta_3\in(P_1\cap P_3)\setminus\{0\}$ and express it as a sum
of the other two nonzero elements of $P_1$. Otherwise choose
$\delta_3\in P_3\setminus(P_1\cap P_2)$; it has a decomposition
with nonzero summands in $P_1$ and $P_2$.
The three corresponding pairs are disjoint and have total sum zero,
contrary to $0\notin\Sigma_6(D)$.

We also need to restrict the sizes of the two possible larger
intersections. Any set of more than four elements in an $H$-coset
has every nonzero element of $H$ as a pair sum. For six elements
$X$ in one coset, complementation therefore gives
\[
 \Sigma_2(X)=H\setminus\{0\},
 \qquad
 \Sigma_4(X)=H\setminus\{\sigma(X)\}.
\]
Its four-term sums meet the pair sums of any other intersection of
size at least three. Consequently, an intersection of size six
excludes any other of size at least three.
For seven elements, the six-element subsets have different sums,
so their four-term sumsets together cover $H$. An intersection of
size seven, or eight, therefore excludes every other intersection
of size at least two.
Finally, two intersections of size five are impossible: the
four-term sums of one meet the pair sums of the other.

An intersection of size $s\ge2$ contributes $s-2$ to
$|D|-\nu_H(D)$. The restrictions just proved leave a total of at
most five, except that a full eight-element coset can contribute six.
This proves both bounds.
\end{proof}

\begin{lemma}\label{rd:mixed}
Suppose that $A$ is independent and $3\le a\le8$. Set
\[
 R_a=2(34-a)-\binom{a-3}{2},
 \qquad \sigma_H=a_1+a_2+a_3.
\]
If $h_H(B),h_H(C)\le2$, and $t,u$ are the numbers of one-element
intersections of $B,C$, respectively, then
\begin{equation}\label{rd:capacity}
 \frac{b(c+u)}2
 \le \frac{bc+bu+ct-tu}{2}\le R_a.
\end{equation}
If also $bc>R_a$, then
\begin{equation}\label{rd:direction}
 \dim\langle\Sigma_2(C)\rangle\le\frac{c+u}{2}\le\frac{R_a}{b},
\end{equation}
and the analogous assertion holds with $B,C$ interchanged.
If $h_H(B)\ge3$, then $h_H(C)=1$ and
\begin{equation}\label{rd:large}
 c\nu_H(B)\le R_a.
\end{equation}
Without assuming that $A$ is independent, the same conclusions
hold with $2m$ in place of $R_a$ whenever $B+C$ meets at most
$m$ cosets of $H$.
\end{lemma}
\begin{proof}
Every element of $B+C$ has a unique representation. Otherwise
two different representations would give a zero-sum subsequence
with two terms from $S_1$ and two from $S_{-1}$.
Consider the array $(b+c)_{b\in B,c\in C}$. Two entries in the same
$H$-coset and in different rows and columns must differ by $\sigma_H$.
Every other difference in $H$ is a sum of at most two of
$a_1,a_2,a_3$, and would complete those four terms to a short
zero-sum subsequence.

This has two consequences. First, a $3\times2$ rectangle within
one quotient coset is impossible: comparing its first entry with
two entries in the other column would make those entries equal.
Thus $h_H(B)\ge3$ forces $h_H(C)=1$.
Second, if either $B_q$ or $C_r$ has size at least two, the sumset
$B_q+C_r$ is the only contribution to the quotient coset $q+r$.
An entry from another contribution would have to differ by
$\sigma_H$ from two distinct entries. If all contributions are
from singleton intersections, at most two of them can occupy one
quotient coset.

Assign a pair of nonempty intersections weight one if both are
singletons, and weight two otherwise. The preceding paragraph
shows that each quotient coset accommodates total weight at most
two. Independence of $A$ gives more precise restrictions.
The images of $a_4,\ldots,a_a$ are independent in $V/H$, and
$B+C$ avoids $\Sigma_{\le 4}^{0}(A)$. The resulting bounds are
\[
\begin{array}{c|c|c}
\text{quotient coset}&\text{number of such cosets}&\text{allowed weight}\\\hline
0\text{ or }a_i+H\ (i\ge4)&a-2&0\\
a_i+a_j+H\ (4\le i<j\le a)&\binom{a-3}{2}&1\\
\text{all remaining cosets}&34-a-\binom{a-3}{2}&2
\end{array}
\]
In the first row, every element uses at most four elements of $A$.
In the second row, all elements except
$a_i+a_j+\sigma_H$ use at most four. In particular, that coset
cannot contain a contribution of weight two. Summing the last
column with these multiplicities gives $R_a$.

If $b=2r+t$ and $c=2s+u$, counting the four types of pairs of
intersections gives total weight
\[
 2rs+2ru+2ts+tu=\frac{bc+bu+ct-tu}{2}.
\]
Subtracting $b(c+u)/2$ leaves $t(c-u)/2\ge0$, proving
\eqref{rd:capacity}. If $bc>R_a$, neither set can have only
singleton intersections. Fix a two-element intersection of $B$
with pair sum $d$. For every two-element intersection of $C$,
with pair sum $e$, the array condition gives $d+e=\sigma_H$.
All these pair sums $e$ are therefore equal.
There are $(c+u)/2$ nonempty intersections of $C$. Their
representatives have difference space of dimension at most
$(c+u)/2-1$; adjoining the common pair sum proves
\eqref{rd:direction}.

If $h_H(B)\ge3$, we have already proved $h_H(C)=1$.
Each element of $C$ contributes total weight $\nu_H(B)$,
which proves \eqref{rd:large}. The argument with total capacity
$2m$ uses only the universal weight bound two per coset, so it
does not require independence of $A$.
\end{proof}

\begin{lemma}\label{rd:change}
Suppose $\dim\langle A\rangle\ge4$, and let $D,E$ be $B,C$ in
either order, with $|D|=N$ and $|E|=M$. Suppose that, for every
coordinate subspace $H$, the condition $h_H(D)\ge3$ implies
$M\nu_H(D)\le R$, where $R$ is fixed. If $(N-5)M>R$, then
$h_H(D)\le2$ for every coordinate subspace $H$.

In particular, when $A$ is independent, this conclusion holds if
\begin{equation}\label{rd:change-condition}
 (N-5)(27-a-N)>R_a.
\end{equation}
For $4\le a\le8$, it holds throughout
$10\le N\le\min\{16,22-a\}$.
\end{lemma}
\begin{proof}
Suppose $h_H(D)\ge3$. The strict inequality gives
$\nu_H(D)\le N-6$, so Lemma~\ref{rd:coset-count} forces $D$
to contain a full coset $F=x+H$.
Choose $a_4\in A\setminus H$, and replace
$H=\langle a_1,a_2,a_3\rangle$ by
$H'=\langle a_1,a_2,a_4\rangle$.
The set $F$ meets exactly two $H'$-cosets, in four elements each.
Both cosets lie in the affine four-space $F+H'$.

There is no third $H'$-coset containing three elements of $D$.
Consequently, every element of $D$ outside $F+H'$ contributes
one to $\nu_{H'}(D)$, while the two cosets inside it contribute
four in total. Since $h_{H'}(D)\ge3$, the same strict inequality
gives $\nu_{H'}(D)\le N-6$. Therefore
\[
 |D\cap(F+H')|=N-\bigl(\nu_{H'}(D)-4\bigr)\ge10.
\]
Choose ten of these elements, forming a set $X$. Lemmas~\ref{rr:cosets}
and~\ref{rr:counting} give
$|\Sigma_3(X)|=r(T_X)\ge F_{10}(T_X)\ge g(10)=20$.
But all these sums belong to the same affine four-space, which
has only $16$ elements. This contradiction proves the first assertion.

For the stated range of $N$, the left side of
\eqref{rd:change-condition} minus $R_a$ is a concave quadratic
in $N$. At the two endpoints its values are
\[
\begin{array}{c|ccccc}
a&4&5&6&7&8\\\hline
N=10&5&3&2&2&3\\
N=\min\{16,22-a\}&17&9&2&2&3
\end{array}
\]
and hence it is positive throughout the interval.
\end{proof}

\begin{lemma}\label{rd:sixteen}
If $A$ is independent and $a=7$ or $8$, then $b,c\le15$.
\end{lemma}
\begin{proof}
By symmetry, suppose $b=16$. Proposition~\ref{rr:mixed} gives
$F_{16}(T_B)\le128$, so Lemma~\ref{rr:counting} leaves
$T_B=d_{16}$ or $e_8\oplus e_8$.
If $T_B=e_8\oplus e_8$, or if the eight pair indices for $d_{16}$
are affinely dependent, the block sumset formulas give a hyperplane
$K\subset V$ such that
\[
 \Sigma_3(B)=V\setminus K,
 \qquad |K\setminus\Sigma_4(B)|\le1.
\]
In the dependent-index case, the only possible even relation among
the eight indices involves all eight; the same formulas apply.
Since $A$ avoids $\Sigma_3(B)$, we have $A\subset K$, and thus
$a=7$ and $c=4$. The set $C$ avoids $\Sigma_4(B)$, so it has at most
one element in $K$. It has at most two outside $K$: three such
elements have sum in $V\setminus K=\Sigma_3(B)$, giving a zero-sum
subsequence of length six. This contradicts $c=4$.

We may therefore assume that $B$ consists of eight parallel pairs
with common pair sum $d\ne0$, and that their indices modulo
$\langle d\rangle$ are affinely independent. In particular, $B$
affinely spans $V$, and $B+d=B$.

We first show that $d$ is a sum of at most two elements of $A$.
The following table gives the intersection sizes of two Hamming
balls whose centers have distance $t\ge3$:
\[
\begin{array}{c|ccc}
\text{dimension, radius}&t=3,4&t=5,6&t\ge7\\\hline
7\text{ and }2&6&0&0\\
8\text{ and }3&38&20&0
\end{array}
\]
To check the table, translate one center to zero. If the other
center has weight $t$, a point with $i$ ones in those $t$
coordinates and $j$ elsewhere belongs to both radius-$r$ balls
exactly when $i+j\le r$ and $t-i+j\le r$.
For example, at distance three the counts are
$3+3=6$ and $2+6(1+5)=38$, respectively.

The set $B+C$ is invariant under translation by $d$ and avoids
$\Sigma_{\le 4}^{0}(A)$. If $a=8$, complementation with respect to the
basis $A$ sends its $bc=48$ elements into the radius-three ball.
Translation invariance puts all of them in the intersection of two
such balls whose centers differ by $d$. Since $48>38$, the table
gives $\wt_A(d)\le2$.

If $a=7$, put $W=\langle A\rangle$. Were $d\notin W$, the set
$B+C$ would avoid both $\Sigma_{\le 4}^{0}(A)$ and its translate by $d$.
These are disjoint sets of size
$1+7+\binom72+\binom73+\binom74=99$.
Only $256-198=58$ elements would remain for the $bc=64$ distinct
mixed sums, which is impossible. Thus $d\in W$.
Since $B$ affinely spans $V$, it meets both $W$-cosets; both
intersections have positive even size because $B+d=B$.
Each element of $C$ consequently gives at least two mixed sums in
$W$, all distinct. Their number is at least $2c=8$.
Complementation in the basis $A$ puts them in two radius-two balls
in $W$, again with centers differing by $d$. Since $8>6$,
we obtain $\wt_A(d)\le2$ in this case as well.

Set $J=\langle\Sigma_2(C)\rangle$, so $\dim J\le c-1$.
We can choose a coordinate subspace $H$ containing $d$ such that
\begin{equation}\label{rd:choose-H}
 \sigma_H+d\notin J.
\end{equation}
Indeed, if $d\in A$, the candidates for $\sigma_H+d$ are the pair
sums of the remaining elements of $A$; their span has dimension
$a-2$. If $d$ is a sum of two elements of $A$, the candidates are
the remaining individual elements, again spanning dimension $a-2$.
In either case $a-2>c-1$, so not all candidates lie in $J$.

Now $h_H(C)=1$. Otherwise two elements of $C$ in the same
$H$-coset would have pair sum $u\in H$, and the mixed-array condition, applied
with a $d$-pair of $B$, would give $u+d=\sigma_H$, contrary to
\eqref{rd:choose-H}.
Every nonempty $B_q$ has even size. None is full, since four pair
indices would then lie in an affine plane of $V/\langle d\rangle$
and would be affinely dependent.
Lemma~\ref{rd:coset-count} therefore gives $\nu_H(B)\ge11$;
this number is even, so it is at least $12$.
There are at least six nonempty intersections of $B$, each of size
at least two. Pairing them with the $c$ singleton intersections of
$C$ requires at least $6c$ distinct quotient cosets, each with
allowed weight two. The table in Lemma~\ref{rd:mixed} provides only
\[
 34-a-\binom{a-3}{2}
 =\begin{cases}21,&a=7,\\16,&a=8.
 \end{cases}
\]
But $6c$ is $24$ or $18$, respectively. This final contradiction
excludes $b=16$.
\end{proof}

\begin{lemma}\label{rd:dense}
Suppose $h_H(C)=1$ and $c\ge13$. Then every $B_q$ has size at
most four, at most two have size at least two, and at most one has
size four. In particular,
\[
 b\le\nu_H(B)+3,
 \qquad \nu_H(B)\le3\ \Longrightarrow\ b\le5.
\]
\end{lemma}
\begin{proof}
We first prove an elementary fact: eleven elements of an affine
four-space have two disjoint triples of each prescribed sum in
that space. Translate the prescribed sum to zero and discard zero
if it occurs. At least ten elements remain. A sum-free subset of
$\F_2^4$ has size at most eight: it is disjoint from its translate
by any of its elements. Hence there is a triple $E$ with
$\sigma(E)=0$.

The plane $P=\langle E\rangle$ has $P\setminus\{0\}=E$.
At least seven of the remaining elements lie outside $P$.
They occupy its three other cosets, each of size four, so two cosets
contain at least three and at least two elements, respectively.
The three elements in the first coset give all three possible
nonzero pair sums in $P$. We can therefore select one pair from
each coset with different sums $u,v\in E$.
Adjoining $u$ to the first pair and $v$ to the second gives two
disjoint zero-sum triples. Translating back proves the fact.

Let $\bar C$ be the image of $C$ in $V/H$. Its at least twelve
nonzero elements form a sum-free set. Indeed, a zero-sum triple
in the quotient would lift to three elements of $C$ whose sum
lies in $H$; at most three elements of $A$ would complete them
to a short zero-sum subsequence.
The Davydov--Tombak result in \cite[Fact~2]{GL10} puts these
nonzero elements in a nonzero coset $L$ of a hyperplane of $V/H$.
The set $L$ is an affine four-space.

Every $q\in L$ has two disjoint three-term representations from
$\bar C\cap L$, by the preceding fact. Every $q\in L+L$ has two
four-term representations sharing exactly one term: fix
$p\in\bar C\cap L$, delete it, and apply the same fact to the at
least eleven remaining elements with target $q+p\in L$.
In each case the two corresponding sums in $V$ are distinct.
Equality would give six distinct elements of $C$ with sum zero.

For five elements $X$ of $B$ in one $H$-coset, complementation of pairs
gives $|\Sigma_3(X)|=7$. If its quotient coset lies in $L$, these
seven sums must meet the two distinct corresponding sums from
$\Sigma_3(C)$, which is impossible.
For four elements $X$ of $B$ in one $H$-coset, put
$U=\{0,a_1,a_2,a_3\}$. We have $|X+U|\ge7$.
Indeed, two missing elements would force $X$ to avoid two
translates of $U$. Two different translates of this simplex
have union of size at least six, leaving room for at most two
elements of $X$, not four.
If the quotient coset of $X$ lies outside $L$, it lies in $L+L$.
The two distinct four-term sums from $C$ in that coset must then
meet $X+U$. This contradicts
\[
 \Sigma_3(B)\cap\Sigma_3(C)=\varnothing,
 \qquad
 \Sigma_4(C)\cap(B+U)=\varnothing,
\]
which follow from the absence of short zero-sum subsequences.
Thus all intersections of size four lie in $L$, and none has
size five or more.

Finally, suppose $B_q,B_r$ are distinct intersections of size at
least two and $q+r\in L+L$. The translates
$q+(\bar C\cap L)$ and $r+(\bar C\cap L)$ each have at least
twelve elements in the same sixteen-element coset, so they meet.
This contradicts the exclusive-contribution property in
Lemma~\ref{rd:mixed}. There is therefore at most one nonsingleton
intersection in each of the two cosets of $L+L$.
There are at most two altogether, and at most one of size four.
They contribute at most $(4-2)+(3-2)=3$ beyond $\nu_H(B)$.
If $\nu_H(B)\le3$ and all intersections are singletons, then
$b=\nu_H(B)\le3$. Otherwise there is one nonsingleton intersection
and at most one further element; hence $b\le4+1=5$.
\end{proof}

\begin{lemma}\label{rd:transverse}
Let $W=\langle A\rangle$.
If $\dim J\le\dim W-3$, some coordinate subspace $H$ satisfies
$H\cap J=0$. Moreover, if $H\cap H'=0$ and
$\nu_H(D),\nu_{H'}(D)\le m$, where $m=3$ or $4$, then
\begin{equation}\label{rd:two-quotients}
 |D|\le2m-3.
\end{equation}
\end{lemma}
\begin{proof}
The image of $W$ modulo $J$ has dimension at least three.
Choose three elements of $A$ with independent images; their span
is the required $H$.

For the second assertion, an $H$-coset and an $H'$-coset meet in
at most one element. If either quotient has all intersections of
size at most two, then $|D|\le m\le2m-3$.
Otherwise consider the $H$-cosets. If two intersections have size
at least two, they already contribute four to $\nu_H(D)$, so
$m=4$ and there are no others. Every $H'$-intersection then has
size at most two, reducing to the preceding case.

The only remaining possibility is one intersection of size at
least three and at most $m-2$ other elements. The large intersection
meets each $H'$-coset at most once. The other elements can create
at most $m-3$ elements in total beyond the first two in each $H'$-coset:
for $m=3$ there is only one other element, and for $m=4$ there
are at most two. Thus
$|D|\le\nu_{H'}(D)+(m-3)\le2m-3$.
\end{proof}

\begin{proof}[Proof of \textup{(P)}]
Interchange $B,C$ if necessary so that $b\ge c$.

\emph{First suppose that $A$ is dependent.}
Every nonempty relation among its elements has weight at least
seven. Two independent relations, together with their sum, would
have total weight at least $21$, whereas every coordinate occurs
zero or twice, giving total weight at most $2a\le16$.
Thus there is a unique nonzero relation, of weight seven or eight.
Put $W=\langle A\rangle$.

If $a=8$, replacing a subset by its symmetric difference with
this relation shows that every element of $W$ is a sum of at most
four elements of $A$. Hence $W=\Sigma_{\le 4}^{0}(A)$, and $B+C$ lies
in the other $W$-coset.
If $a=7$, then $\dim W=6$ and $W=\Sigma_{\le 3}^{0}(A)$.
The images of $B,C$ in $V/W$ are disjoint, because $B+C$ avoids
$W$. Each image also misses a nonzero element of $V/W$:
representatives of all three nonzero elements would form a triple
whose sum lies in $W$, and at most three elements of $A$ would
complete a short zero-sum subsequence.

In this latter case, a class meeting $W$ has at most two elements
there and at most one in each other coset. Three elements in $W$,
or one in $W$ and two in another coset, would give the same
forbidden type of triple. Since $b,c\ge4$, such a class must have
intersection sizes $2,1,1$. The other class would consist of
sixteen elements in the sole remaining coset. Its triple sums
would lie in a $64$-element affine space, contrary to
$r(T)\ge g(16)=112$. Thus neither class meets $W$.
Their disjoint supports among the three nonzero quotient elements
then put $B+C$ in at most two $W$-cosets.

For $a=8$, one $W$-coset contains sixteen cosets of a coordinate
subspace. For $a=7$, two $W$-cosets contain $2\cdot8=16$ such
cosets. In both cases Lemma~\ref{rd:mixed} therefore applies
with total capacity $32$ for every coordinate subspace.
The permitted sizes give $(b-5)c\ge33$ (the minimum is $44$
when $a=7$ and $33$ when $a=8$). Lemma~\ref{rd:change} gives
$h_H(B)\le2$ for every $H$.

Suppose $h_H(C)\le2$ for some $H$. Since $bc\ge48>32$,
the dimension bound gives
$\dim\langle\Sigma_2(C)\rangle\le\lfloor32/b\rfloor\le3$.
As $\dim W\ge6$, Lemma~\ref{rd:transverse} supplies a coordinate
subspace transverse to this difference space. Its quotient map
is injective on $C$, so the mixed bound becomes $bc\le32$, a
contradiction. We conclude that $h_H(C)\ge3$ for every $H$.
The mixed bound now gives $\nu_H(C)\le\lfloor32/b\rfloor\le3$.

Choose two coordinate subspaces $H,H'$ with $H\cap H'=0$,
using six independent elements of $A$. Relative to $H$, the set
$C$ lies in one intersection of size at least three and at most
one additional element. The large intersection meets each
$H'$-coset at most once; the additional element can only create
an intersection of size two. This contradicts $h_{H'}(C)\ge3$.

\emph{Now suppose that $A$ is independent.}
Lemma~\ref{rd:sixteen} excludes $b=16$ for $a\ge7$.
For $a=8$, the possibility $b=15,c=4$ is also excluded by packing:
\[
 256\ge 1+8+\binom82+\binom83+15\cdot4+g(15)
 =93+60+112=265.
\]
It follows that $c\ge5$. For $a\ge4$, every class of size at
least ten falls within the range of Lemma~\ref{rd:change}; its
intersection with each coordinate coset therefore has size at most two.

If $a=3$, the weighted bound cannot hold with
$h_H(B),h_H(C)\le2$. Indeed, $b+c=24$ and $8\le b,c\le16$.
For even $b,c$, its minimum is $bc/2\ge64$; for odd $b,c$,
its minimum is $(bc+23)/2\ge79$. Both exceed $R_3=62$.
Interchange $B,C$ so that $h_H(B)\ge3$; then $h_H(C)=1$.
Combining $b\le\nu_H(B)+6$ with $c\nu_H(B)\le62$ gives
\[
 c(18-c)=c(b-6)\le62.
\]
Since $8\le c\le16$, this forces $14\le c\le16$.
Lemma~\ref{rd:dense} improves $b\le\nu_H(B)+6$ to
$b\le\nu_H(B)+3$, giving the contradiction
$62\ge c(b-3)=c(21-c)\ge80$.

We may now assume $4\le a\le8$. Since $b\ge10$, we have
$h_H(B)\le2$ for every coordinate subspace.
Suppose first that $h_H(C)\le2$ for at least one $H$, and set
$J=\langle\Sigma_2(C)\rangle$. The permitted sizes give
$bc\ge70>R_a$, so
Lemma~\ref{rd:mixed} gives $\dim J\le R_a/b$.
If $\dim J\le a-3$, a coordinate subspace transverse to $J$
would have $h_H(C)=1$, and hence $bc\le R_a$, again impossible.
Consequently
\begin{equation}\label{rd:dimension-test}
 b(a-2)\le R_a.
\end{equation}

There is also a useful parity restriction. In
\eqref{rd:capacity}, $t\equiv b\pmod2$ and $u\equiv c\pmod2$.
The expression is nondecreasing in each of $t,u$, so its minimum is
\[
 P(b,c)=2\lceil b/2\rceil\lceil c/2\rceil
 -\begin{cases}1,&b,c\text{ both odd},\\0,&\text{otherwise}.
 \end{cases}
\]
The conditions $P(b,c)\le R_a$ and \eqref{rd:dimension-test}
leave only two possibilities:
\[
\begin{array}{c|c|l}
a&R_a&\text{consequence}\\\hline
4&60&b+c=23\text{ gives }P(b,c)\ge64\\
5&57&(b,c)=(14,8)\text{ or }(16,6)\\
6&53&b\le13\text{ gives }P(b,c)\ge56\\
7&48&b\ge10\text{ gives }b(a-2)\ge50\\
8&42&b\ge10\text{ gives }b(a-2)\ge60
\end{array}
\]
For completeness, at $a=5$ the pairs with $b=11,\ldots,16$
have parity costs $71,60,69,56,63,48$, respectively.

In either surviving case, $h_H(C)\ge3$ is impossible for any $H$.
It would force $h_H(B)=1$. For $(b,c)=(14,8)$, we would then
have $\nu_H(C)\le\lfloor57/14\rfloor=4$ and, by
Lemma~\ref{rd:dense}, $c\le7$. For $(b,c)=(16,6)$, we would
have $\nu_H(C)\le3$ and the same lemma would give $c\le5$.
Thus both classes have intersections of size at most two for every $H$.

For $(b,c)=(14,8)$, inequality \eqref{rd:capacity} reads
$14u+8t-tu\le2$. Since $t,u$ are even, it forces $t=u=0$.
For $(b,c)=(16,6)$, it reads $16u+6t-tu\le18$, forcing $u=0$.
Hence every nonempty intersection of the $14$-element class in
the first case, or of the $6$-element class in the second, has
size two. Their common pair sum $d_H$ is nonzero and satisfies
$D+d_H=D$ for that class $D$.
Such a period is unique: two distinct nonzero periods would
generate a four-element translation group, making $|D|$ divisible
by four, whereas $|D|$ is $14$ or $6$.
The same nonzero period must therefore belong to every coordinate
subspace. Their intersection is zero, since $A$ is independent
and $a=5$. This contradiction excludes the possibility
$h_H(C)\le2$ for any $H$.

We have proved that $h_H(C)\ge3$ for \emph{every} coordinate
subspace. We finish by applying the same bounds twice.
Lemma~\ref{rd:change} gives $c\le9$, so $b\ge18-a$.
For $a=4,\ldots,8$, one has $5(18-a)>R_a$; hence
the mixed bound gives $\nu_H(C)\le4$ for every $H$.
If $a=4,5$, then $b\ge13$ and $h_H(B)=1$, so
Lemma~\ref{rd:dense}, with the classes interchanged, gives $c\le7$.
If $a\ge6$, two disjoint coordinate subspaces and
Lemma~\ref{rd:transverse} give $c\le5$.

These improvements give $b\ge20-a$ in the first case and
$b\ge22-a$ in the second. In either case $4b>R_a$, so now
$\nu_H(C)\le3$ for every $H$.
Repeat the preceding bounds. For $a=4,5$,
Lemma~\ref{rd:dense} gives $c\le5$, whereas $b\le16$ requires
$c=27-a-b\ge11-a\ge6$.
For $a\ge6$, Lemma~\ref{rd:transverse} gives $c\le3$, whereas
we already know $c\ge5$. These final contradictions prove
\textup{(P)}.
\end{proof}

\section*{Acknowledgements}

This work was supported by the National Natural Science Foundation
of China (grant nos. 12371335 and 12271520).

\end{document}